# On Nelson Goodman's final formulae for Primary Complexity

Godofredo Iommi Amunátegui

*Instituto de Física, Pontificia Universidad Católica de Valparaíso, Casilla 4059,*

*Valparaíso, Chile*

giommi@ucv.cl

In the *Structure of Appearance* and in *Problems and Projects*, Nelson Goodman has constructed a theory of complexity whose elements are the predicates of a system. One of his main results is a closed formula to evaluate $v\ n-pl$, the maximum complexity value of an n-place predicate. Up to this day no procedure has been found to reduce $v\ n-pl$, i.e., to write it in terms of (n-1)-place predicates. In this article we propose elementary combinatorial algorithms to carry out such a reduction.

# I

In these pages we wish to consider anew the pioneering analysis of Nelson Goodman on complexity [Goodman, 1966, pp. 63-123; 1972, pp.33-40 and pp.275-354]. Our study concerns a central theoretical issue which is cleary displayed in the following two quotations from Quine and Goodman, respectively:

(i) "Yet it is possible in general to subject general terms to striking formal condensation. If the assumed universe of objects includes at least a modest fund of classes –actually none of more than two members are required for the purpose– then it can be shown that any vocabulary (finite or infinite) of general terms (absolute on relative) is reducible by paraphrase to a single dyadic term." [Quine, 1960, p.232] [1]

(ii) "Many people have supposed that my method of measuring the simplicity of a set of extralogical predicates is undermined by such results as the reduction Quine cites here. If any set of predicates can be replaced by a simple two-place predicate, does this not show that no set of predicates is more complex than a single two-place predicates?" [Goodman, 1972, p.319]

Goodman's assurance is based on Lars Svenonius'proof that it is not in general possible to replace an n-place predicate of individuals by any set of predicates of individuals having fewer than n places each [Svenonius, 1955]. After, grosso modo, more than fifty years the concept of complexity has been extended to new areas of knowledge and its very meaning has been, somewhat, modified. The literature on the subject is, by now, rather large, to say the least[2]. Goodman's theory has been criticized and generalized, however its basic tenet, i.e., the irreducibility of the formulae for primary complexity, remains to this day unquestioned by any subsequent consideration.

---
[1] See also Goodman [1972, p.319]

The aim of this work is to show, by means of elementary combinatorial arguments, that Goodman's robust results are reducible, that is, the complexity value of an n-predicate turns out to be reducible to the complexity value of an expression with (n-1)-place predicates, and, moreover, that a reducibility chain can be established in such a manner that an n-place predicate may be written in terms of 1-place predicates.

# II

Complexity and simplicity are twin concepts. To clarify the point the following lines are helpful: "Our present concern is with the logical or structural simplicity of bases which pertains directly to the degree of systematization of theories founded upon them. In summary, then, we want to find a way of measuring the structural simplicity of the set of undefined extralogical terms of a theory or a system. That is, we want to be able to assign to any such set of terms a number that will indicate the complexity of that set accordingly, one significant aspect of the complexity of the theory" [Goodman, 1972, p. 284]. In this section, Goodman's ideas will be outlined, almost *ad litteram*, in order

---

[2] Let us mention three articles where Goodman's contribution is considered explicitly: Wimsatt [1972], Turney [1989] and, specially, Richmond [1996].

to fix their scope and meaning.

1. All scientific activity amounts to the invention of and the choice among systems of hypotheses. A system is achieved just to the extent that the basic vocabulary and set of first principles used in dealing with the given subject matter are simplified. Systematization is the same thing as simplification of basis [Goodman, 1972, p. 279].

2. A theory is a system of statements. Hence the simplicity of the set of concepts, or the vocabulary of terms, employed in these statements must be taken into account. Some words and symbols like "and", "or", "if ... then", "all", "some", "=" or translations of these, are -properly- logical terms. The extralogical basis of a system consists of all its primitives that are not in the list of basic "logical terms". To adopt a term as primitive is to introduce it into a system without defining it [Goodman, 1972, p.283; 1966, p.63].

3. The formal simplicity of bases is related to those differences among predicates that are expressible by using logical terms in addition to the predicates themselves. Among the extralogical terms of a system may be property-terms like "... is green", and relation-terms of various degrees like "... is larger than ... " and "... lies between ... and ... ". The examples given are respectively a one-place, a two-place and a three-place predicate. Only extralogical predicates are counted as contributing to the complexity -value of a basis. The complexity-value of the basis b is denoted by

"vb" [Goodman, 1972, p. 283; 1966, p.68].

4. The complexity-value of every extralogical predicate is a positive integer and the complexity-value of a basis is the sum of the values of the extralogical predicates in it. A basis consisting only of logical predicates has the value 0. Concerning notation, use of square brackets will result in a name for the relevant specification abbreviated within them. Thus the full reading of $[2-pl.irref.; two1-pl]$ is "the relevant specification basis consisting of one 2-place irreflexive predicate and two 1-place predicates". All extralogical 1-place predicates have equal value. Goodman sets the value of extralogical 1-place predicates at 1: $v[1-pl.] = 1$. The complexity of any basis can be proved to be equal to that of a certain basis consisting only of thoroughly irreflexive predicates ("... is a parent of ..." is irreflexive) [Goodman, 1972, p.290; 1966, pp. 73-75].

5. The number $h_j$ of $n-j-pl.$ irreflexive predicates [Goodman, 1966, p. 83] and the complexity value of any basis [Goodman, 1966, p.104] were obtained by Goodman with the help, respectively, of P. Savage and N.J. Fine. Hence these results will be called through this article, Goodman-Savage numbers and Goodman-Fine numbers: The general schema for an n-place predicate is: $v[n-pl] = v[n-pl.irref.; h_1 n-1pl.irref.; \ldots h_{n-2} 2-pl.irref.]+1$ where for each j, the number $h_j$ of n-j-pl. irreflexive predicates is

$$\sum_{r=0}^{n-j}(-1)^r \frac{(n-j-r)^n}{(n-j-r)!\,r!} \quad \text{(Goodman-Savage numbers)} \tag{1}$$

Besides, the complexity-value of any basis whatever is given by the general formula:

$$v(n-pl) = \sum_{k=1}^{n}\sum_{r=0}^{k} \frac{(-1)^r (2k-1)(k-r)^n}{r!\,(k-r)!} \quad \text{(Goodman-Fine numbers)} \tag{2}$$

These formulae cover irreflexive predicates only [Goodman, 1966, pp. 82-83]. The complexity-values considered dealt with so far are to be considered *primary* complexity-values. The maximum values for bases may be thought of as measure of potential complexity. The *secondary* complexity-value of a basis is the difference between the complexity-potential and the primary complexity of the basis [Goodman, 1966, p.105].

# III

Svenonius considers the definition of relations in terms of other relations and studies how a definition of a relation $R_1$ in terms of another relation $R_2$ works in a given universe U. In particular, he deals with the question [Svenonius, 1955, p. 235]: "Is every relation of the form $s_1$ definable in terms of some relation of the forms $s_2$?" Definability $D(s_1,s_2)$ is a key concept. For our purpose, Svenonius results may be exposed, in a rather cursory way, as follows:

(a) If $s_1$ and $s_2$ are relation forms, and for some value of the argument the function $M_{s_1}(e)$ is greater than the function $M_{s_2}(e)$, then $D(s_1,s_2)$ is not valid [Svenonius, 1955, p. 249].

(b) If $U = {1,2,3,...}$, $s_1 = s(3\text{-place relation})$ and $s_2 = s(2\text{-place relation})$, then $D(s_1,s_2)$ is not valid [Svenonius, 1955, p.249]. Some details of this last theorem's proof are relevant: $s_1$ contains the five forms (xxx), (xxy), (xyx), (yxx), (xyz), and $M_{s_2}(e) = 2^5$, to $s_2$ correspond the two forms (xx), (xy), and $M_{s_2}(e) = 2^2$. Since $2^5 > 2^2$, $D(s_1,s_2)$ is not valid.

(c) If $k$ and $n$ are integers, $m=n-1$, $s_1 = s$ (n-place relation) and $s_2 = s$ (sequence of $k$ m-place relations), then $D(s_1,s_2)$ is not valid [Svenonius, 1955, p. 249].

Svenonius states that his developments are related to the criterion used by Goodman for degrees of complexity [Svenonius, 1955, p. 235]. The most important consequence of such an interpretation is that if $D(s_1,s_2)$ is not valid then the basis $s_1$ cannot *be reduced* to a basis $s_2$. Goodman says that Svenonius "has proved that it is not in general possible even to replace an n-place predicate of individuals by any set of predicates of individuals having fewer than n places each" [Goodman, 1972, p. 319]. This is a direct reply to Quine.

It seems interesting, at this stage, to note that Goodman is aware of a theoretical bias appearing in Svenounius' approach to complexity: "So to admit all structure-classes as relevant kinds is to obtain a measure of power rather of what I have called complexity" [Goodman, 1972, p. 312].

# IV

The Goodman-Savage numbers $h_{n,n-j}$ (where $j = 0, 1, ..., n-1$) may be arranged in a triangular matrix. Instead of labeling its entries as $(n, n-j)$ we shall denote them by $(n, n+j)$ with the same values for $j$. Hence we have:

$$h_{n,j+1} = \sum_{r=0}^{j+1} -1^r \frac{j+1-r^n}{j+1-r\ !r!}$$

multiplying this equation by $1 = \frac{j+1\ !}{j+1\ !}$

$$h_{n,j+1} = \sum_{r=0}^{j+1} -1^r \frac{j+1-r^n}{j+1-r\ !r!} \cdot \frac{j+1\ !}{j+1\ !} = \frac{1}{j+1\ !}\sum_{r=0}^{j+1} \frac{(j+1)!}{j+1-r)!r!}\ j+1-r^n$$

So

$$h_{n,j+1} = \frac{1}{j+1\ !}\sum_{r=0}^{j+1} -1^r\ C_r^{j+1}\ j+1-r^n$$

Let $j+1 = m$ and $k = m-r$; moreover $C_r^{j+1} = C_{j+1-r}^{j+1}$, then

$$h_{n,m} = \frac{1}{m!}\sum_{k=0}^{m} (-1)^{m-k} C_k^m k^n \qquad (3)$$

That is: $h_{n,m} = S_n^{(m)} \equiv$ Stirling numbers of the second kind[3].

**Proposition 1:** The Goodman-Savage numbers are identical with the Stirling numbers of the second kind and therefore correspond to the number of ways of partitioning a set of n elements into m non-empty subsets.

Table I: The Goodman-Savage numbers

$\underrightarrow{j+1}$  $j = 0, 1, \ldots, n-1$

$n\downarrow$

| | | | | | | |
|---|---|---|---|---|---|---|
| 1 | | | | | | |
| 1 | 1 | | | | | |
| 1 | 3 | 1 | | | | |
| 1 | 7 | 6 | 1 | | | |
| 1 | 15 | 25 | 10 | 1 | | |
| 1 | 31 | 90 | 65 | 15 | 1 | |
| 1 | 63 | 301 | 350 | 140 | 21 | 1 |

---

[3] See, for instance, Abramowitz and Stegun [1972, pp. 824-825].

Where $j = 0,1, \ldots, 6$ and $n = 1,2, \ldots, 7$

This identity between $h_{n,m}$ and $S_n^{(m)}$ allows us to use the following relation of this kind of Stirling numbers. Namely:

$$(j+1)h_{n,j+1} = h_{n+1,j+1} - h_{n,j} \qquad (4)$$

Equation (2) gives a procedure to evaluate the Goodman-Fine numbers which may also be ordered as a lower-left triangular matrix:

Table II: The Goodman-Fine numbers

$\underrightarrow{j+1}$

| | | | | | | |
|---|---|---|---|---|---|---|
| 1 | | | | | | |
| 1 | 3 | | | | | |
| 1 | 9 | 5 | | | | |
| 1 | 21 | 30 | 7 | | | |
| 1 | 45 | 125 | 70 | 9 | | |
| 1 | 93 | 450 | 455 | 135 | 11 | |
| 1 | 189 | 1505 | 2450 | 1260 | 231 | 13 |

$n\downarrow$

Analyzing Equation (2) we remark that, for a given *n*, the sum over *k* indicates the maximum value for a k-place predicate which can be written as $v_k(n-pl)$.

For instance, for *n* = 5 and *k* = 3:

$$v(5-pl) = 5\sum_{r=0}^{3}(-1)^r \frac{(3-r)^5}{(3-r)!r!} = 125$$

In general:

$$v(n-pl) = \sum_{k=1}^{n} v_k(n-pl) \quad (5)$$

where

$$v_k(n-pl) = (2k-1)\sum_{r=0}^{k} \frac{(-1)^r (k-r)^n}{r!(k-r)!} \quad (6)$$

As an example:

$$v(4-pl) = v_1(4-pl) + v_2(4-pl) + v_3(4-pl) + v_4(4-pl) = 1+21+30+7 = 59$$

Equation (2) may be expressed as:

$$v(n-pl) = \sum_{j=0}^{n-1}\sum_{r=0}^{j+1} (-1)^r (2j+1) \frac{(j+1-r)^n}{(j+1-r)!r!}$$

But

$$h_{n,j+1} = \sum_{r=0}^{j+1} (-1)^r \frac{(j+1-r)^n}{(j+1-r)!r!}$$

And we obtain

$$v(n-pl) = \sum_{j=0}^{n-1} (2j+1) h_{n,j+1} \qquad (7)$$

Equation (7) states that

**Proposition 2:** The Goodman-Fine numbers may be written in terms of the Goodman-Savage numbers.

Example:

$$v(5-pl) = \sum_{j=0}^{4} (2j+1) h_{5,j+1}$$

$$v(5-pl) = 1 \cdot h_{5,1} + 3h_{5,2} + 5h_{5,3} + 7h_{5,4} + 9h_{5,5}$$

$$v(5-pl) = 1 \cdot 1 + 3 \cdot 15 + \cdot 5 \cdot 25 + 7 \cdot 10 + 9 \cdot 1 = 1 + 45 + 125 + 70 + 9 = 250$$

The case $k = n$ of (6) is interesting:

$$v_n(n-pl) = (2n-1)\sum_{r=0}^{n}(-1)^r \frac{(n-r)^n}{(n-r)!r!}$$

Let us consider the sum:

$$\sum_{r=0}^{n}(-1)^r \frac{(n-r)^n}{(n-r)!r!} \cdot \frac{n!}{n!} = \sum_{r=0}^{n}(-1)^r \frac{n!}{(n-r)!r!} \cdot \frac{(n-r)^n}{n!} = \sum_{r=0}^{n}(-1)^r C_{n-r}^{n} \frac{(n-r)^n}{n!}$$

But $n! = \sum_{r=0}^{n}(-1)^r C_{n-r}^{n}(n-r)^n$

$$\therefore \frac{\sum_{r=0}^{n}(-1)^r C_{n-r}^{n}(n-r)^n}{n!} = 1$$

Finally:

a) $v_n(n-pl) = 2n-1$. In such a manner we recover Goodman's formula 3.741 [Goodman, 1966, p. 103].

b) From (6) and (7) we derive:

$$v_k(n-pl) = (2k-1)h_{n,k} \qquad (8)$$

Besides by means of a combination of eq. (4) and eq. (7), we get:

**Proposition 3:** $v(n-pl)$ may be written in terms of $h_{n-1,j}$, i.e.,

$$v(n-pl) = \sum_{j=1}^{n-1} A_j h_{n-1,j} \qquad (9)$$

Where $A_j = 2j^2 + j + 1$ and $n-1 \geq j$

Ex.: $v(7-pl) = 4h_{6,1} + 11h_{6,2} + 22h_{6,3} + 37h_{6,4} + 56h_{6,5} + 79h_{6,6}$.

Equation (9) establishes a first type of reduction: symbolically $n \to n-1$. Such a reduction procedure may be applied to each of the steps of the chain $n \to n-1 \to n-2 \to \ldots \to 2 \to 1$. Each one of these reduction cases requires a longer-but straightforward-calculation. We display here some characteristic results:

$$v(n-pl) = \sum_{j=1}^{n-2} (jA_j + A_{j+1}) h_{n-2,j} \qquad n-2 \geq j$$

$$v(n-pl) = \sum_{j=1}^{n-3} (j^2 A_j + (2j+1) A_{j+1} + A_{j+2}) h_{n-3,j} \qquad n-3 \geq j$$

$$\vdots \qquad \qquad \vdots$$

$$v(n-pl) = \left( \sum_{j=1}^{n-1} A_j h_{n-1,j} \right) h_{1,1}$$

For Goodman all extralogical 1-place predicates have equal value and his third postulate $v(n-pl) = 1$ sets the value of extralogical 1-place predicates [Goodman, 1966, p.75]. From the recurrence relation (4): $1 \cdot h_{1,1} = h_{2,1} - 0 = 1$, so $v(n-pl) = h_{1,1}$, and $v(n-pl)$ may be expressed as:

$$v(n-pl) = K v(1-pl) \text{ where } K = \sum_{j=1}^{n-1} A_j h_{n-1,j} \text{ is an integer.} \qquad (10)$$

For example, for n=4

$$v(4-pl) = 59 \cdot v(1-pl)$$

**Proposition 4:** The maximum complexity value of any n-place predicate breaks down into K one-piece predicates, where K is an integer.

Using the same combinatorial tools, we can display the maximum complexity value of a n-place predicate in terms of an n-1 place predicate:

**Proposition 5:** The complexity value $v(n-pl)$ may be written as

$$v(n-pl) = 3v[(n-1)-pl] + \sum_{j=1}^{n-1} C_j h_{n-1,j} \text{ where } j=1,2,3,\ldots (n-1) \text{ and}$$

$$C_j = 1 + (j-1)(1+2(j-2)). \tag{11}$$

For instance:

$$v(6-pl) = 3(5-pl) + 1 \cdot h_{5,1} + 2h_{5,2} + 7h_{5,3} + 16h_{5,4} + 29h_{5,5}$$

From Proposition 5 it is possible to deduce another expression for $v(n-pl)$ in terms of $v(1-pl)$. It suffices to put, successively, $n-1, n-2, n-3, \ldots, 1$ into the left-hand side of equation (11). We have then:

**Proposition 6:** $v(n-pl) = 3^{n-1} v(1-pl) + R$

Where $R = 3^{n-2}C_1h_{1,1} + 3^{n-3}C_1h_{2,1} + 3^{n-3}C_2h_{2,2} + ... + 3^0C_1h_{n-1,1} + 3^0C_2h_{n-1,2} + ... + 3^0C_{n-1}h_{n-1,n-1}$ (12)

Example:

$v\ 5-pl = 3^4 v\ 1-pl + 3^3C_1h_{1,1} 3^3C_1h_{2,1} + 3^2C_2h_{2,2} + 3^1C_1h_{3,1} + 3^1C_2h_{3,2} + 3^1C_3h_{3,3} + 3^0C_1h_{4,1} + 3^0C_2h_{4,2} + 3^0C_3h_{4,3} + 3^0C_4h_{4,4}$

# V

Svenonius[4] [1955, p. 299] has proved that if $s_1 = s$ (3-place relation) and $s_2 = s$ (2-place relation), the Definability $D\ s_1,s_2$ is not valid. Hence, not every basis of the kind $s_1$ is always replaceable by some basis $s_2$. Proposition 5 indicates that $v\ 3-pl = 3v\ 2-pl + C_1h_{2,1} + C_2h_{2,2} = 15$ (see Table II). In general for $s_1 = s$ ($n-1$ -place relation) and $s_2 = s$ $n-1$ -place relation), even if $D\ s_1,s_2$ is not valid, $v\ n-pl$ may be reducible to $v\lfloor n-1-pl \rfloor$ plus an integer. It appears that Goodman's Complexity and Svenonius' Definability are non-equivalent concepts.

At the beginning of this article we treat with some detail Goodman's theoretical certitude concerning the irreducibility of the complexity value of an n-place predicate. Let us quote other statements dealing with such a topic:

---

[4] See Section III (b) of the present article. In Svenonius' proof, the five forms of $s_1$ and the two forms of $s_2$ correspond, respectively, to the sum of the Goodman-Savage numbers $h_{3,1} + h_{3,2} + h_{3,3}$ and to $h_{2,1} + h_{2,2}$ (see Table I).

"Various writers have devised ways of showing that every basis can in some sense "be reduced to" some basis consisting of a single 2-place predicate. Actually none of these devices establishes the replaceability required by my way of measuring complexity. But the threat remained that some new invention might do the trick or something like it" [Goodman, 1972, p. 312].

According to his view, Svenonius has proved "that it is not the case that every basis can always be replaced by some basis consisting of a single 2-place predicate, or even by some basis consisting of any finite numbers of 2-place and 1-place predicates" [Goodman, 1972, p.312].

As a matter of fact, no "new invention" is needed "to do the trick". Propositions 4,5 and 6 allow us to write the maximum complexity value of an n-place basis in terms of the complexity values of r-place predicates ( r = 1,2,3, …, n-1)[5]

---

[5] The main point of our combinatorial argument is contained in equations (3) and (4).

# Appendix

If we denote the matrices corresponding, respectively, to the Goodman-Fine numbers and to the Goodman-Savage numbers by $M_{G-F}$ and $M_{G-S}$, for a given **k**, we have $M_{G-S} \cdot \Delta = M_{G-F}$, where $\Delta$ is a diagonal matrix whose elements are 1,3,5,2$n$-1. For instance, for $n = 4$

| 1 |   |   |   |   | 1 |   |   |   |   |   | 1 |    |    |   |
|---|---|---|---|---|---|---|---|---|---|---|---|----|----|---|
| 1 | 1 |   |   |   |   | 3 |   |   |   |   | 1 | 3  |    |   |
| 1 | 3 | 1 |   |   |   |   | 5 |   |   | = | 1 | 9  | 5  |   |
| 1 | 7 | 6 | 1 |   |   |   |   | 7 |   |   | 1 | 21 | 30 | 7 |

Note that $M_{G-S}\Delta = M_{G-F}$ is nothing by Equation (7) written in matrix form.

# Aknowledgments

I am grateful to Marcus Rossberg for having brought to my attention Quine´s review of the "Structure of Appearance" (The Journal of Philosophy, Vol. 48, No. 18 (1951), pp. 556-563) and to Bruno Merello Encina for helpful comments.

This work was supported by Fondo Nacional de Desarrollo Científico y Tecnológico [1120019].